\documentclass[12pt]{article}
\usepackage{amsmath,latexsym}
\usepackage{amscd}
\usepackage{amsfonts}
\usepackage{mathrsfs}
\usepackage{stmaryrd,bbm}
\usepackage{dsfont}
\usepackage{graphics}
\usepackage{graphicx}
\usepackage{amssymb}
\newcommand{\mapform}[5]{\begin{array}{ccrcl}
#1 & : & #2 & \longrightarrow & #3 \\
& & #4 & \longmapsto & #5
\end{array}}
\DeclareMathSymbol{\mP}{\mathbin}{AMSb}{'120}
\newcommand{\qed}{\hfill $\blacksquare$}

\newtheorem{example}{Example}

\newtheorem{theorem}{Theorem}
\newtheorem{corollary}{Corollary}
\newtheorem{lemma}{Lemma}
\newtheorem{definition}{Definition}

\begin{document}
\title{Geometric Modeling and Regularization of Algebraic Problems}
\author{Zhonggang Zeng \thanks{Department of Mathematics,
Northeastern Illinois University, Chicago, Illinois 60625, USA.
~~email:~{\tt zzeng@neiu.edu}. 
~Research is supported in part by NSF under grant DMS-1620337.}}

\newcommand{\CC}{\mathbbm{C}}
\newcommand{\cC}{\mathcal{C}}
\newcommand{\cK}{\mathcal{K}}
\newcommand{\cP}{\mathcal{P}}
\newcommand{\cF}{\mathcal{F}}
\newcommand{\cU}{\mathcal{U}}
\newcommand{\cV}{\mathcal{V}}
\newcommand{\cW}{\mathcal{W}}
\newcommand{\cM}{\mathcal{M}}
\newcommand{\cN}{\mathcal{N}}
\newcommand{\bdu}{\mathbf{u}}
\newcommand{\bdv}{\mathbf{v}}
\newcommand{\bdw}{\mathbf{w}}
\newcommand{\bdx}{\mathbf{x}}
\newcommand{\bdy}{\mathbf{y}}
\newcommand{\bdz}{\mathbf{z}}
\newcommand{\bda}{\mathbf{a}}
\newcommand{\bdb}{\mathbf{b}}
\newcommand{\bdc}{\mathbf{c}}
\newcommand{\bdd}{\mathbf{d}}
\newcommand{\bde}{\mathbf{e}}
\newcommand{\bdf}{\mathbf{f}}
\newcommand{\bdg}{\mathbf{g}}
\newcommand{\bdh}{\mathbf{h}}
\newcommand{\bdi}{\mathbf{i}}
\newcommand{\bdj}{\mathbf{j}}
\newcommand{\bdm}{\mathbf{m}}
\newcommand{\bdn}{\mathbf{n}}
\newcommand{\bdk}{\mathbf{k}}
\newcommand{\bdl}{\mathbf{l}}
\newcommand{\bdr}{\mathbf{r}}
\newcommand{\bdp}{\mathbf{p}}
\newcommand{\bdq}{\mathbf{q}}
\newcommand{\bds}{\mathbf{s}}
\newcommand{\bdo}{\mathbf{0}}
\newcommand{\bdF}{\mathbf{F}}
\newcommand{\h}{{{\mbox{\tiny $\mathsf{H}$}}}}
\newcommand{\dm}{\dim}
\newcommand{\dg}{\deg}
\newcommand{\codim}{\mathrm{codim}}
\newcommand{\la}{\lambda}
\newcommand{\al}{\alpha}
\newcommand{\bt}{\beta}
\newcommand{\gm}{\gamma}
\newcommand{\dl}{\delta}
\newcommand{\Dl}{\Delta}
\newcommand{\eps}{\varepsilon}
\newcommand{\sg}{\sigma}
\newcommand{\rank}[1]{\mathrm{rank}\left(\,#1\,\right)}
\newcommand{\ranka}[2]{\mathrm{rank}_{#1}\left(\,#2\,\right)}
\newcommand{\GCD}{\mathrm{gcd}}
\newcommand{\dist}[1]{\mathrm{dist}\left(\,#1\,\right)}

\maketitle

\begin{abstract}
Discontinuity with respect to data perturbations is common in algebraic
computation where solutions are often highly sensitive. 
Such problems can be modeled as solving systems of equations at given
data parameters.
By appending auxiliary equations, the models can be formulated to 
satisfy four easily verifiable conditions so that the data form complex
analytic manifolds on which the solutions maintain their structures and
the Lipschitz continuity. 
When such a problem is given with empirical data, solving the system
becomes a least squares problem whose solution uniquely exists and enjoys 
Lipschitz continuity as long as the data point is in a tubular neighborhood 
of the manifold.
As a result, the singular problem is regularized as a well-posed
computational problem.
\end{abstract}

\vspace{-4mm}
\section{Introduction}

\vspace{-4mm}
Computational problems with extremely high sensitivities beyond finite 
bounds are known to be {\em ill-posed}. 
Such problems are abundant in algebraic computation and also referred to
as being {\em singular}.
Some of the most basic algebraic problems are ill-posed, such as matrix 
ranks and subspaces, solutions of singular linear systems, polynomial 
greatest common divisors and factorizations, defective eigenvalues and 
Jordan Canonical Forms.
Those are the problems we inevitably encounter in symbolic, numeric and 
hybrid computation.
Based on the current state of knowledge, however, it is inaccurately believed 
by many that such problems are impossible to solve from empirical data or 
using floating point arithmetic.
Pessimistic outlooks are abundant in the literature (emphasis added): 
``The moral is to {\em avoid floating point solutions} of singular 
systems''\cite[page 218]{Meyer}.
``The difficulty is that the JCF {\em cannot be computed} using floating point 
arithmetic.
A single rounding error may cause some multiple eigenvalue to become 
distinct or vice versa, altering the entire structure''
\cite{MolVanL}.
``A {\em dramatic deterioration of the accuracy} must therefore be
expected''\cite[page 300]{ueberhuber}.
``[S]mall variations in the [data]
will result in large variations in the [solution].
{\em There is no hope of computing such an object} in a stable way''%
\cite[page. 128]{HorJoh}.
``[Although such an object] is of fundamental theoretical 
importance {\em it is of little use} in practical computations, being generally 
very difficult to compute''\cite[page 52]{bar-cam}.
``[So] that [it] is little used in numerical applications''%
\cite[page. 128]{HorJoh}.

Are the solutions of those problems really sensitive to data perturbations as 
alleged?
In a legendary technical report \cite{kahan72}, Kahan argues that it is a 
``misconception'' to consider multiple roots of polynomials hypersensitive,
points out that polynomials and matrices form heuristic ``pejorative 
manifolds'' preserving root multiplicities and Jordan structures respectively, 
and proves that the sensitivities of roots and eigenvalues are 
bounded if the perturbation is constrained to preserve the multiplicity.
This insight opens a possible pathway for accurate solution of such singular 
problems.

In this paper, we establish conditions for modeling an algebraic problem
as a nonlinear system of equations in the form of solving 
$\bdf(\bdu,\bdv)=\bdo$ for the variable $\bdv$ at a fixed data value 
$\bdu$ so that we can rigorously verify that the data form a complex analytic 
manifold on which the solution maintains a certain algebraic structure and 
enjoys Lipschitz continuity.

The data of a hypersensitive problem forming smooth manifolds is 
crucial in the analysis and regularization the problems since its solution 
is of bounded sensitivity with respect to data on the manifold. 
We further extend this inherent stability beyond the manifold into its
{\em tubular neighborhood}.
When the problem data are given as empirical, we have a data point 
near the manifold in the data space.
Assuming the data are reasonably accurate so that the point remains in the 
tubular neighborhood, the Tubular Neighborhood Theorem established in this
paper ensures the projection from the data point to the manifold uniquely 
exists and enjoys Lipschitz continuity. 
Consequently, the singular problem can be regularized as a well-posed
least squares problem that is accurately solvable from empirical data.

The geometric modeling and regularization from this perspective lead to robust 
algorithms such as those in accurate computation of multiple roots 
\cite{zeng-mr-05}, greatest common divisors \cite{KalYanZhi06,ZengGCD}, 
polynomial factorizations \cite{WuZeng,ZengAIF}, defective eigenvalue 
problems \cite{pseudoeig} and singular linear systems \cite{Zeng2019}.
These algorithms are implemented in our software package 
{\sc NAClab} \cite{naclab}.

Geometric theories and methods have been applied in algebraic computing in 
many works such as 
\cite{AbsMahSep,cor-gal,Ded96,dem-edel, edelman-elmroth-kagstrom_1,
edelman-elmroth-kagstrom_2, EdeAriSmi,NumAlgGeo}.
However, the tremendous advantage of tubular neighborhoods has not yet 
been harnessed partly because a general tubular neighborhood theorem for 
complex analytic manifolds is apparently unavailable in the literature of 
differential topology.
Specifically taylored to the application of solving ill-posed algebraic
problems in this paper, we prove a weak but sufficient version of the
tubular neighborhood theorem for complex analytic manifolds in Euclidean 
spaces using the techniques of nonlinear least squares.
The theorem and the proof fills a gap in the regularization theory of solving 
ill-posed algebraic problems and complete the works of numerical 
factorization \cite{WuZeng,ZengAIF} and numerical greatest common divisors of 
polynomials \cite{zeng-dayton,ZengGCD}. 

\vspace{-4mm}
\section{Preliminaries}

\vspace{-4mm}
The space of $n$-dimensional vectors of complex numbers is denoted by $\CC^n$ 
with the Euclidean norm $\|\cdot\|_2$.
General vector spaces are denoted by, say $\cV$, $\cW$ in which 
vectors are denoted by boldface lower case letters while ~$\bdo$ ~is a zero 
vector.
Any norm $\|\bdv\|$ is understood as the specified norm in 
the space where $\bdv$ belongs.

For a holomorphic mapping $F:\Omega\subset\CC^n\rightarrow\CC^m$, we may 
designate a variable name, say $\bdz$, for $F$ and denote $F$ as
$\bdz\mapsto F(\bdz)$.
The Jacobian matrix of $F$ at any $\bdz_0\in\Omega$ is denoted by 
$F_\bdz(\bdz_0)$.
Let $\cV$ and $\cW$ be vector spaces with isomorphisms 
$\phi:\cV\rightarrow\CC^n$ and $\psi:\cW\rightarrow\CC^m$. 
Assume $\bdg$ is a mapping from an open subset $\Sigma$ of $\cV$ to $\cW$
with a representation $\bdz\mapsto G(\bdz)$ where 
$G:\phi(\Sigma)\subset\CC^n\rightarrow \CC^m$ such that 
$\bdg=\psi^{-1}\circ G\circ \phi$.
We say $\bdg$ is holomorphic in $\Sigma$ if its representation $G$ is 
holomorphic in $\phi(\Sigma)$. 
Denoting the variable of $\bdg$ as, say $\bdv$, the {\em Jacobian} of $\bdg$ 
at any particular $\bdv_0\in\Sigma$ is defined as the linear transformation 
$\bdg_\bdv(\bdv_0):\cV\rightarrow\cW$ in the form of 
\[ \bdv\longmapsto\bdg_\bdv(\bdv_0)(\bdv) :=
\psi^{-1}\big( G_\bdz(\phi(\bdv_0))\,\phi(\bdv)\big).
\]
The Jacobian $\bdg_\bdv(\bdv_0)$ as a linear transformation is invariant 
under change of bases.
Let $G_\bdz(\bdz_0)^\h$ and $G_\bdz(\bdz_0)^\dagger$ be the Hermitian 
transpose and the Moore-Penrose inverse of the Jacobian matrix 
$G_\bdz(\bdz_0)$ respectively where $\bdz_0=\phi(\bdv_0)$.
If we further assume the isomorphisms $\phi$ and $\psi$ are isometric, 
namely $\|\phi(\bdv)\|_2 =\|\bdv\|$ and $\|\psi(\bdw)\|_2 = \|\bdw\|$ 
for all $\bdv\in\cV$ and $\bdw\in\cW$, then $\bdg_\bdv(\bdv_0)^\h$ and 
$\bdg_\bdv(\bdv_0)^\dagger$ are well-defined as 
\[ \bdg_\bdv(\bdv_0)^\h = \phi^{-1}\circ G_\bdz(\bdz_0)^\h\circ\psi 
~~~~\mbox{and}~~~
\bdg_\bdv(\bdv_0)^\dagger  = \phi^{-1}\circ G_\bdz(\bdz_0)^\dagger\circ\psi
\]
that are invariant under isometric isomorphisms.
A mapping $\bdf$ is holomorphic in a non-open domain $\Pi\subset\cV$ if 
there is an open subset $\Omega$ of $\cV$ containing $\Pi$ and a holomorphic 
mapping $\bdg$ defined in $\Omega$ such that $\bdf(\bdz)\equiv\bdg(\bdz)$ for 
all $\bdz\in\Pi$.
For a holomorphic mapping $(\bdu,\bdv)\mapsto\bdf(\bdu,\bdv)$, its Jacobian
at $(\bdu_0,\bdv_0)$ is denoted by $\bdf_{\bdu\bdv}(\bdu_0,\bdv_0)$ and its
partial Jacobian with respect to, say $\bdv$, is denoted by 
$\bdf_\bdv(\bdu_0,\bdv_0)$.

\vspace{-4mm}
\section{Complex analytic manifolds}

\vspace{-4mm}
%
For our applications, 
we consider complex analytic manifolds in normed vector spaces in the 
following definition.

\begin{definition}[Complex Analytic Manifold] \label{d:cam}
Let $\cU$ be a finite-dimensional normed vector space over $\CC$. 
A subset $\Pi$ of $\cU$ is a {\em complex analytic manifold} of dimension $m$ 
if there is an $m$-dimensional normed vector space $\cV$ over $\CC$ and, for 
every $\bdu \in \Pi$, there is an open neighborhood $\Sigma$ of $\bdu$ in 
$\cU$ and a holomorphic mapping $\phi$ from $\Sigma \cap \Pi$ onto an open 
subset $\Lambda$ of $\cV$ with a holomorphic inverse. 
The dimension deficit $\dm(\cU)-m$ is called the {\em codimension} of $\Pi$ 
in $\cU$ denoted by $\codim(\Pi)$.
\end{definition}

The term {\em manifold} in this paper refers to a complex analytic manifold 
in the sense of Definition~\ref{d:cam}.
As we shall elaborate in case studies in \S\ref{s:mod}, algebraic problems 
whose solutions possess certain algebraic structures can often be modeled as 
a system of nonlinear equations in the form of {\em solving 
$\bdf(\bdu,\bdv) = \bdo$ ~for the variable $\bdv$ ~at the given data 
parameter value $\bdu$}. 
The following theorem establishes four basic conditions for such a model
so that the data points form a manifold. 
The theorem simplifies the tedious process of establishing a manifold to 
verifying the four conditions.

\vspace{-4mm}
\begin{theorem}[Geometric Modeling Theorem]\label{t:cam}
A subset $\Pi$ is a complex analytic manifold in a 
normed vector space $\cU$ over $\CC$ if and only if there are normed
vector spaces $\cV$ and $\cW$ over $\CC$ with
\[ \dm(\cV) ~\le~ \dm(\cW) ~\le~ \dm(\cU)+\dm(\cV) ~<~ \infty
\]
such that, at every $\bdu_0\subset\Pi$, there is a holomorphic 
mapping $(\bdu,\bdv)\mapsto\bdf(\bdu,\bdv)$ from an open domain
$\Omega\subset\cU\times\cV$ to $\cW$ with the properties below:

\vspace{-6mm}
\begin{itemize}\parskip0mm
\item[\em (i)] There is a $\bdv_0\in\cV$ such that 
$\bdf(\bdu_0,\bdv_0)=\bdo$.
\item[\em (ii)] $\bdf_{\bdu\bdv}(\bdu_0,\bdv_0)$ is 
surjective and $\bdf_{\bdv}(\bdu_0,\bdv_0)$ is injective.
\item[\em (iii)] $\bdf(\bdu,\bdv)=\bdo$ implies $\bdu\in\Pi$.
\item[\em (iv)]  For every open neighborhood $\Dl$ of $\bdv_0$ in $\cV$, there 
is an open neighborhood $\Sigma$ of $\bdu_0$ in $\cU$ such that every 
$\bdu\in\Sigma\cap\Pi$ corresponds to a unique $\bdv\in\Dl$ with 
$\bdf(\bdu,\bdv)=\bdo$. 
\end{itemize}

\vspace{-6mm}
Under these conditions, we have $\codim(\Pi) = \dm(\cW)-\dm(\cV)$.
\end{theorem}

\vspace{-4mm}
{\bf Proof.} 
Let $\Pi$ be a manifold in $\cU$ with $\cW\,=\,\cU$ as in 
Definition~\ref{d:cam}, the mapping 
\[
(\bdu,\bdv)\mapsto\bdf(\bdu,\bdv) ~=~ \bdu-\phi^{-1}(\bdv)
\]
from 
$\Sigma\times\Lambda$ in $\cU\times\cV$ to $\cW$ satisfies conditions (i)-(iv).

Conversely, assume $\bdf$ satisfies all the specified conditions
and we proceed to prove $\Pi$ is a manifold in $\cU$.
From property (ii), we can write
\[ \cU ~=~\hat\cU\oplus\check\cU ~~~~\mbox{with}~~~
\dm(\check\cU)+\dm(\cV) ~=~ \dm(\cW),
\] 
regard $\cU$ as $\hat\cU\times\check\cU$ and consider $\bdf$
as $(\hat\bdu,\check\bdu,\bdv)\mapsto\bdf(\hat\bdu+\check\bdu,\bdv)$
from the domain $\Omega$ in $\hat\cU\times\check\cU\times\cV$ to $\cW$ 
so that $\bdf_{\check\bdu\bdv}(\hat\bdu_0,\check\bdu_0,\bdv_0)$ is invertible
where $\hat\bdu_0+\check\bdu_0=\bdu_0$.
By the Implicit Mapping Theorem \cite{Tay00}, there is a neighborhood 
$\Lambda\times\Dl$ of $\big(\hat\bdu_0,\,(\check\bdu_0,\bdv_0)\big)$ in 
$\hat\cU\times(\check\cU\times\cV)$, holomorphic mappings 
$\bdg:\Lambda\subset\hat\cU\rightarrow\check\cU$ and
$\bdh:\Lambda\subset\hat\cU\rightarrow\cV$ such that 
$(\check\bdu_0,\, \bdv_0)= (\bdg(\hat\bdu_0),\,\bdh(\hat\bdu_0))$ and
$\bdf(\hat\bdu+\check\bdu,\bdv\big)=\bdo$ for 
$\big(\hat\bdu,(\check\bdu,\bdv)\big)$ in $\Lambda\times\Dl$ 
if and only if $(\check\bdu,\bdv)=(\bdg(\hat\bdu),\,\bdh(\hat\bdu))$.
Without loss of generality, we assume 
$\Lambda\times\Dl=\Omega$ since we can redefine $\bdf$ with 
a restricted domain.

Let $\psi$ be the holomorphic mapping
$\hat\bdu\mapsto(\hat\bdu,\,\bdg(\hat\bdu))$ from $\Lambda\subset
\hat\cU$ to $\hat\cU\times\check\cU$.
Then $\psi(\Lambda)\subset\Pi$ since 
$(\psi(\hat\bdu),\bdh(\hat\bdu))=(\hat\bdu,\bdg(\hat\bdu),\bdh(\hat\bdu))$ 
is in $\bdf^{-1}(\bdo)$ for all $\hat\bdu\in\Lambda$.
We also have $\psi(\hat\bdu_0)=(\hat\bdu_0,\bdg(\hat\bdu_0))
=(\hat\bdu_0,\check\bdu_0)$. 
Let $\tilde\Dl=\{\bdv\in\cV \,|\, (\check\bdu_0,\bdv)\in\Dl\}$
that is an open neighborhood $\bdv_0$ in $\cV$.
By the condition (iv), there is an open neighborhood $\Sigma$ of 
$(\hat\bdu_0,\check\bdu_0)$ in $\hat\cU\times\check\cU$ such that every 
$(\hat\bdu,\check\bdu)\in\Sigma\cap\Pi$ corresponds to a unique 
$\bdv\in\tilde\Dl$ with $\bdf(\hat\bdu,\check\bdu,\bdv)=\bdo$.
Denote $\tilde\Lambda=\psi^{-1}(\Sigma)$ that is open in $\hat\cU$ and define 
the holomorphic mapping $\phi:(\hat\bdu,\check\bdu)\mapsto\hat\bdu$
from $\Sigma\subset\hat\cU\times\check\cU$ to $\hat\cU$.
Clearly $\phi\circ\psi(\hat\bdu)=\phi(\hat\bdu,\bdg(\hat\bdu)) = \hat\bdu$ 
for all $\hat\bdu\in\tilde\Lambda\subset\Lambda$.
Furthermore, for every $(\hat\bdu,\check\bdu)\in\Sigma\cap\Pi$, there is a 
unique $\bdv\in\tilde\Dl$ with $\bdf(\hat\bdu,\check\bdu,\bdv)=\bdo$ so that 
$(\check\bdu,\bdv)=(\bdg(\hat\bdu),\bdh(\hat\bdu))$.
Namely 
\[ \psi\circ\phi(\hat\bdu,\check\bdu) ~=~ \psi(\hat\bdu) ~=~
(\hat\bdu,\bdg(\hat\bdu)) ~=~ (\hat\bdu,\check\bdu).
\]
Consequently, the subset $\Pi$ is a manifold in $\cU=\hat\cU\times\check\cU$ 
of dimension $\dm(\hat\cU)$ that equals to $\dm(\cU)+\dm(\cV)-\dm(\cW)$. 
\qed

Assuming the model of solving $\bdf(\bdu,\bdv)=\bdo$ for $\bdv$ at the given
data $\bdu$ is properly formulated so that the conditions of 
Theorem~\ref{t:cam} are satisfied, the solution $\bdv$ is locally Lipschitz 
continuous with respect to the data $\bdu$ on the manifold.

\vspace{-4mm}
\begin{corollary}\label{c:mlip}
Using the notations in Theorem~\ref{t:cam}, assume $\bdf$ satisfies the 
condition {\em (i)-(iv)}.
Further assume $\cU$, $\cV$ and $\cW$ are normed and the isomorphisms from 
$\cV$ and $\cW$ to $\CC^{\dm(\cV)}$ and $\CC^{\dm(\cW)}$, respectively, 
are isometric. 
Then there is an open neighborhood $\Omega_0$ of $\bdu_0$ in $\cU$ such that, 
for every fixed parameter $\bdu_1\,\in\,\Omega_0\,\cap\,\Pi$, the equation 
$\bdf(\bdu_1,\bdv)\,=\,\bdo$ has a unique solution $\bdv_1\,\in\,\cV$ and
\begin{equation}\label{v1v0}
  \|\bdv_1-\bdv_0\| \le \|\bdf_\bdv(\bdu_0,\bdv_0)^\dagger\|\,
\|\bdf_\bdu(\bdu_0,\bdv_0)\|\,\|\bdu_1-\bdu_0\| + o(\|\bdu_1-\bdu_0\|).
\end{equation}
\end{corollary}

\vspace{-4mm}
{\em Proof.}
Using the notations in the proof of Theorem~\ref{t:cam}, we have
$\bdf(\hat\bdu+\bdg(\hat\bdu),\bdh(\hat\bdu))\,\equiv\,\bdo$ for 
$\hat\bdu\,\in\,\Lambda$, implying the linear transformation
\[ \bdf_{\hat\bdu}(\hat\bdu_0+\check\bdu_0,\bdv_0)+
\bdf_{\check\bdu}(\hat\bdu_0+\check\bdu_0,\bdv_0)\circ
\bdg_{\hat\bdu}(\hat\bdu_0)+
\bdf_{\bdv}(\hat\bdu_0+\check\bdu_0,\bdv_0)\circ\bdh_{\hat\bdu}(\hat\bdu_0)
\]
maps $\hat\bdu_1-\hat\bdu_0$ to $\bdo$ from all $\hat\bdu_1\in\hat\cU$.
Furthermore, from 
\begin{align*}
\bdv_1-\bdv_0 &\,=\,\bdh_{\hat\bdu}(\hat\bdu_0)\,(\hat\bdu_1-\hat\bdu_0)+h.o.t
\mbox{~~~and} \\
\check\bdu_1-\check\bdu_0 &\,=\, 
\bdg_{\hat\bdu}(\hat\bdu_0)\,(\hat\bdu_1-\hat\bdu_0)+h.o.t.
\end{align*}
where $h.o.t.$ denotes the sum of higher order terms of 
$\hat\bdu_1-\hat\bdu_0$. 
Since $\cV$ and $\cW$ are isometrically isomorphic to $\CC^{\dm(\cV)}$ and 
$\CC^{\dm(\cW)}$ respectively so that $\bdf_\bdv(\bdu_0,\bdv_0)^\dagger$ is 
well-defined and we have (\ref{v1v0}) from
\begin{align*}  \bdv_1-\bdv_0
& =  -\bdf_\bdv(\hat\bdu_0+\check\bdu_0,\bdv_0)^\dagger \,
\big(\bdf_{\hat\bdu}(\hat\bdu_0+\check\bdu_0,\bdv_0)\,(\hat\bdu_1-\hat\bdu_0)
\\
& \;\;\;\;\;\;\;\;\;\;\;\;\;\;\;\;\;\;\;\;\;\;\;\; 
+ \bdf_{\check\bdu}(\hat\bdu_0+\check\bdu_0,\bdv_0)\, 
(\check\bdu_1-\check\bdu_0)\big) + h.o.t \\
&= -\bdf_\bdv(\hat\bdu_0+\check\bdu_0,\bdv_0)^\dagger \,
\bdf_{\hat\bdu \check\bdu}(\hat\bdu_0+\check\bdu_0,\bdv_0)\,
\big((\hat\bdu_1,\check\bdu_1)-(\hat\bdu_0,\check\bdu_0)\big)+h.o.t \\
&= -\bdf_\bdv(\bdu_0,\bdv_0)^\dagger \,
\bdf_\bdu(\bdu_0,\bdv_0)\, (\bdu_1-\bdu_0)+h.o.t  ~~~~~~~~~\mbox{\qed}
\end{align*}

The solution of a singular problem is known to be infinitely sensitive to 
{\em arbitrary} perturbations. 
In \cite{kahan72}, Kahan discovers an inherently bounded stability under 
perturbations
{\em constrained} on certain heuristically conceived ``pejorative manifolds''
for the root-finding and the eigenvalue problems.
Theorem~\ref{t:cam} rigorously establishes the conditions 
for modeling general algebraic problems so that data points indeed form 
manifolds on which the solutions maintain certain structures.
Corollary~\ref{c:mlip} further quantifies the bounded sensitivity 
on those manifolds.
More importantly, the bounded sensitivity can be extended beyond the manifold 
into its tubular neighborhood, making it possible to harness the stability 
in practical computation from empirical data as we shall elaborate in 
\S\ref{s:tnt}.

\vspace{-4mm}
\section{Geometric modeling case studies}\label{s:mod}

\vspace{-4mm}
Algebraic problems are often phrased in a pattern of finding a certain 
solution at a data point, such as ``find the kernel of a matrix'', ``find 
the greatest common divisor of a polynomial pair'', ``find the Jordan 
Canonical Form of a matrix'', ``find the factorization of a polynomial''.
The data point can usually be represented as a vector $\bdu=\hat\bdu$ in a 
vector space $\cU$. 
The key to the geometric analysis and the accurate solution of those
problems is to model the solution as a vector $\bdv$ in a vector 
space $\cV$ in a zero-finding problem: 
\begin{equation}\label{fuv0} 
\mbox{At } \,\hat\bdu\in\cU, \mbox{ solve the equation }
 \,\bdf(\hat\bdu,\bdv) = \bdo \,\mbox{ for } \bdv\,\in\,\cV
\end{equation}
where $\bdf:(\bdu,\bdv)\mapsto\bdf(\bdu,\bdv)$ is a holomorphic 
mapping from an open domain $\Omega$ in $\cU\times\cV$.
By adding proper auxiliary equations, the model can be set up so that the 
mapping $\bdf$ satisfies the conditions (i)-(iv) in Theorem~\ref{t:cam}.
Consequently, a collection of the data points at which the solutions possess 
a specific algebraic structure can be established as a 
{\em structure-preserving manifold}, making it possible to 
apply the Tubular Neighborhood Theorem (Theorem~\ref{t:tnt}).
The model (\ref{fuv0}) also enables computation of an approximate
solution as the least squares solution 
$\bdv\,=\,\tilde\bdv$ of the equation $\bdf(\tilde\bdu,\bdv)\,=\,\bdo$.
We elaborate such geometric modeling in case studies in this section.

\vspace{-4mm}
\subsection{The matrix rank-revealing problem}\label{s:rr}

\vspace{-4mm}
In $\CC^{m\times n}$ of $m\times n$ matrices of complex entries with the 
Frobenius norm $\|\cdot\|_{_F}$, the subset 
\[ \cC^{m\times n}_r ~:=~ \big\{ A\in\CC^{m\times n} \,\big|
\,\rank{A} = r\big\}
\]
is a manifold of codimension $(m-r)(n-r)$.
This result is proved in \cite{dem-edel} and can be easily verified via 
using Theorem~\ref{t:cam} as follows.

Let $O$ and $I$ denote the zero and identity matrices, respectively, 
in $\CC^{m\times n}$.
At a matrix $A\in\CC^{m\times n}$ of rank-$r$, consider the 
rank-revealing problem as finding the kernel $\cK(A)$ of dimension 
$n-r$. 
The fundamental equation is $G\,X=O$ for $X\in\CC^{n\times (n-r)}$
at the data point $G=A$.
The crucial auxiliary equation that ensures proper modeling under 
Theorem~\ref{t:cam} can be derived from the fact that, for almost all 
$C\in\CC^{n\times(n-r)}$, there is an $N\in\CC^{n\times(n-r)}$ whose columns 
form a basis for $\cK(A)$ such that $C^\h\,N=I$.
Finding the kernel of $A$ can then be modeled as a zero-finding problem:
\[  \mbox{Solve\,\, }
\bdf(A,X) = (O,O) \mbox{\,\, for\,\, } X\in\CC^{n\times(n-r)}
\]
where, with a fixed parameter $C\in\CC^{n\times(n-r)}$, the mapping
\begin{equation}\label{rrmap}
\mapform{\bdf}{\Omega\subset\CC^{m\times n}\times\CC^{n\times (n-r)}}{
\CC^{(n-r)\times(n-r)}\times\CC^{m\times(n-r)}}{(G,\,X)}{(C^\h X-I,\,G\,X)}.
\end{equation}
Here $\Omega$ is an open neighborhood of $(A,N)$ and, for every 
$(G,X)\in\Omega$, we have $\|A-G\|_{_F}<\|A^\dagger\|_2^{-1}$.
Clearly $\bdf(A,N)=(O,O)$ and $(G,X)\in\bdf^{-1}(O,O)$ implies $G$ has the 
desired algebraic structure of rank\,$r$, leading to the condition (i) and 
(iii) of Theorem~\ref{t:cam}.
The Jacobian 
\[ \bdf_{G X}(A,N) ~:~ (G,X)\,\longmapsto\, (C^\h\,X,G\,N+A\,X)
\]
is surjective since both $C$ and $N$ are of full rank $n-r$.
The partial Jacobian 
\[ \bdf_{X}(A,N) ~:~ X\,\longmapsto\, (C^\h\,X,A\,X)
\]
is injective since $(C^\h\,X,A\,X)=(O,O)$ implies $X=N\,T$ for a certain
\[ T\,\in\,\CC^{(n-r)\times(n-r)}, ~~~~O ~=~ C^\h\,X~=~T ~~~~\mbox{and}~~~
X ~=~ O, 
\]
leading to the condition (ii) of Theorem~\ref{t:cam}.
Furthermore, every matrix $G\in\cC^{m\times n}_r$ sufficiently close to $A$ 
corresponds to a matrix $X\in\CC^{n\times(n-r)}$ whose column span $\cK(G)$ 
and $C^\h\,X=I$ so $\bdf(G,X)=(O,O)$ and $\|X-N\|_{_F}$ can be as small as 
we wish, validating the condition (iv) of Theorem~\ref{t:cam}. 
As a result, the subset $\cC^{m\times n}_r$ is a manifold of codimension 
\[ \codim(\cC^{m\times n}_r) \,=\,
(n-r)^2+m\,(n-r)-n\,(n-r)\,=\,(m-r)\,(n-r).\]
The subset $\cC_r^{m\times n}$ for every $r$ is a structure-preserving manifold
for the rank-revealing problem and the desired solution (rank and kernel) is
modeled in the vector $X\,\in\,\CC^{n\times (n-r)}$ as the zero of the mapping
$X\,\mapsto\,\bdf(\hat{G},X)$ at $\hat{G}\in\cC_r^{m\times n}$.

\vspace{-4mm}
\subsection{The root-finding problem}\label{s:rf}

\vspace{-4mm}
A polynomial can be considered as a data vector in the vector space $\mP_n$ 
of polynomials with degrees up to $n$ and the norm 
\[\|a_0+a_1\,x+\cdots+a_n\,x^n\|\,:=\, \|(a_0,a_1,\ldots,a_n)\|_2\]
that makes $\mP_n$ isometrically isomorphic to $\CC^{n+1}$.
The complete root-finding problem of a polynomial is equivalent to its 
factorization.
For any positive integers $\ell_1+\cdots+\ell_k=n$, denote
\begin{align}
\cF_{\ell_1\cdots\ell_k} \,:=\, \label{Fell} 
\big\{ \al\,(x-z_1)^{\ell_1}\cdots (x-z_k)^{\ell_k} \big|~
\al,z_1,\ldots,z_k\,\in\,\CC, z_i\ne z_j, ~\forall i\ne j\big\}.
\end{align}
Every polynomial $p\in\mP_n$ belongs to one of such a subset in which the 
factorization structure is preserved.
The root-finding problem of $p$ becomes calculating the distinct 
roots $\bdz_1,\ldots,z_k$ and multiplicities $\ell_1,\ldots,\ell_k$.
At any $p\in\cF_{\ell_1\cdots\ell_k}$ with leading coefficient $u_0$
and distinct roots $u_1,\dots,u_k$ of multiplicities 
$\ell_1,\ldots,\ell_k$ respectively, the root-finding problem of $p$ 
can thus be modeled as identifying $\cF_{\ell_1\cdots\ell_k}$ and solving a 
zero-finding problem in the form of the modified Vie\'te's equation
\begin{equation}\label{phzp0}
\phi(\bdz,p) = 0 \mbox{\,\, for\,\, } \bdz=(z_0,z_1,\ldots,z_k)\in\CC^{k+1}
\end{equation}
with the holomorphic mapping from $\Omega\subset\CC^{k+1}\times\mP_n$ to
$\mP_n$ 
\begin{equation}\label{rfmap}
\phi ~:~ (\bdz,\,g) ~\longmapsto~
z_0\,(x-z_1)^{\ell_1}\cdots (x-z_k)^{\ell_k}-g
\end{equation}
where $\Omega$ is an open neighborhood of $\bdu\,=\,(u_0,u_1,\ldots,u_k)$ 
in $\CC^{k+1}$ in which every $\bdy\,=\,(y_0,y_1,\ldots,y_k)\,\in\,\Omega$ 
implies $(y_0,y_{i_1},\ldots,y_{i_k})$ $\not\in\,\Omega$ whenever the 
permutation $(i_1,\ldots,i_k)\,\ne\,(1,\ldots,k)$.
Such a geometric modeling leads to the geometric insight in the following
theorem along with a proof that is made simple by Theorem~\ref{t:cam}.
The theorem sets the foundation for the accurate solution of root-finding
problem in the presence of multiple roots.
The theorem is proposed in \cite{ZengAIF} by this author with an incomplete 
proof due to necessary abbreviation under the page limit.

\vspace{-4mm}
\begin{theorem}\label{t:rfmf}
The subset $\cF_{\ell_1\cdots\ell_k}$ is a complex analytic manifold in
$\mP_n$ of codimension $n-k$ where $n=\ell_1+\cdots+\ell_k$.
\end{theorem}

\vspace{-4mm}
{\bf Proof.} 
For any $p=u_0\,(x-u_1)^{\ell_1}\ldots(x-u_k)^{\ell_k}\in
\cF_{\ell_1\cdots\ell_k}$ with distinct roots $u_1,\ldots,u_k$, define $\phi$ 
as in (\ref{rfmap}) at $p$ so $\phi(\bdu,p)=0$ where $\bdu=(u_0,\ldots,u_k)$.
For any $g\in\mP_n$, we have $\phi_{\bdz g}(\bdu,p)(\bdo,g)\equiv-g$,
implying $\phi_{\bdz g}(\bdu,p)$ is surjective.
With a proof nearly identical to that of Theorem 3.3 in \cite{zeng-mr-05}, 
the partial Jacobian $\phi_{\bdz}(\bdu,p)$ is injective.
Moreover, the continuity of polynomial roots with respect to the coefficients
ensures the condition (iv) of Theorem~\ref{t:cam} is satisfied, concluding
the proof.  \qed

We call $\cF_{\ell_1\cdots\ell_k}$ a {\em factorization manifold} in $\mP_n$.
Factorization manifolds serve as structure-preserving manifolds for 
polynomials in $\mP_n$.
The desired factorization is represented by the vector $(z_0,z_1,\ldots,z_k)$ 
in $\CC^{k+1}$ in the zero-finding model (\ref{phzp0}).
The root-finding problem is thus equivalent to identifying the factorization
manifold $\cF_{\ell_1\cdots\ell_k}$ along with the zero-finding problem
(\ref{phzp0}).

Modeling the factorization problem for polynomials including multivariate 
cases is given in \cite{WuZeng} where the proof of the Factorization Manifold 
Theorem can be substantially simplified by citing Theorem~\ref{t:cam} rather 
than essentially mirroring its proof.

\vspace{-4mm}
\subsection{The greatest common divisor problem}\label{s:gcd}

\vspace{-4mm}
We say two polynomials are $\sim$-equivalent if they are constant multiples
of each other.
For every $(p,q)\in \mP_m\times\mP_n$, let $\GCD(p,q)$ denote the greatest 
common divisor (GCD) of $p$ and $q$ as an equivalent class under $\sim$. 
The subset $\cP^k_{m,n}$ defined as
\[  \big\{ (p,q)\in\mP_m\times\mP_n
\big| \dg(p)=m,\dg(q)=n,\dg(\GCD(p,q)) = k\big\}
\]
is a manifold of codimension $k$ in $\mP_m\times\mP_n$ where $\dg(\cdot)$ is 
the degree of any polynomial $(\cdot)$, as asserted in \cite{ZengGCD}.
To establish this result, we model the GCD computation as a zero-finding 
problem.
At any particular $(\hat{p},\hat{q})\in\cP^k_{m,n}$, there is a
~$(u,v,w)=(\hat{u},\hat{v},\hat{w})$ satisfying the quations 
$u\,v-p=u\,w-q=0$ at the data $(p,q)=(\hat{p},\hat{q})$ with 
$\hat{u}\in\GCD(\hat{p},\hat{q})$.
To ensure proper modeling, we need an auxiliary equation
$r\odot u=\bt\ne0$ for almost all $r\in\mP_k$ such as a random polynomial 
where $\odot$ is the dot-product between two polynomials defined as the 
dot-product between the corresponding coefficient vectors.
Using such $r$ and $\bt$ as parameters, the GCD problem of the pair
$(\hat{p},\hat{q})$ can be modeled as identifying the GCD degree $k$
and
\[ \mbox{Solve\,\,~} \psi(u,v,w,\hat{p},\hat{q}) = (0,0,0)
\mbox{\,\,~for~\,\,} (u,v,w)\,\in\,\mP_k\times\mP_{m-k}\times\mP_{n-k}
\]
with the holomorphic mapping 
\begin{align}
\mapform{\psi}{
\Omega\subset\mP_k\times\mP_{m-k}\times\mP_{n-k}\times \mP_m\times\mP_n}{
\CC\times\mP_m\times\mP_n}{(u,v,w,p,q)}{(r\odot u-\bt,\,u\,v-p,\,u\,w-q)}
\label{gcdmap}
\end{align}
where $\Omega$ is an open neighborhood of 
$(\hat{u},\hat{v},\hat{w},\hat{p},\hat{q})$ in 
$\mP_k\times\mP_{m-k}\times\mP_{n-k}\times\mP_m\times\mP_n$
such that every $(u,v,w,p,q)\in\Omega$ satisfies $\dg(p)=m$, $\dg(q)=n$ and 
$\dg(u)=k$ with the pair $(v,w)$ being coprime.
Clearly $\psi(\hat{u},\hat{v},\hat{w},\hat{p},\hat{q})=(0,0,0)$.
The Jacobian 
\[
\psi_{uvwpq}(\hat{u},\hat{v},\hat{w},\hat{p},\hat{q}) \,:\,
(u,v,w,p,q)\,\mapsto\,(r\odot u,\, \hat{u}\,v+u\,\hat{v}-p, 
\hat{u}\,w+u\,\hat{w}-q)
\]
can be easily verified to be surjective.
The injectivity of the partial Jacobian 
\[
\psi_{uvw}(\hat{u},\hat{v},\hat{w},\hat{p},\hat{q})\,:\,
(u,v,w)\,\mapsto\,(r\odot u,\, \hat{u}\,v+u\,\hat{v}, 
\hat{u}\,w+u\,\hat{w})
\]
is established by \cite[Corollary 4.1]{ZengGCD}.
At every $(u,v,w,p,q)\in\Omega$, the equality 
\[ \psi(u,v,w,p,q) ~=~ (0,0,0)
\]
implies $(p,q)\in\cP^k_{m,n}$.
It is also a straightforward verification that, for every polynomial pair
$(p,q)\in \cP^k_{m,n}$ sufficiently close to $(\hat{p},\hat{q})$, there is a 
unique GCD triplet
$(u,v,w)\in\mP_k\times\mP_{m-k}\times\mP_{n-k}$ such that 
$\psi(u,v,w,p,q)=(0,0,0)$ with the distance
$\|(u,v,w)-(\hat{u},\hat{v},\hat{w})\|$ as small as we wish.
By Theorem~\ref{t:cam}, the subset $\cP^k_{m,n}$ is a 
manifold in $\mP_m\times\mP_n$ of the codimension
\[ 
\dm(\CC\times\mP_m\times\mP_n)-\dm(\mP_k\times\mP_{m-k}\times\mP_{n-k})
= k.
\]
Each manifold amoung $\cP_{m,n}^0,\,\cP_{m,n}^1,\,\ldots,\,
\cP_{m,n}^{\min\{m,n\}}$ preserves a GCD structure (degree) for polynomial 
pairs on it.

\vspace{-4mm}
\subsection{The Jordan Canonical Form problem}\label{s:jcf}

\vspace{-4mm}
The collection of $n\times n$ matrices with a fixed structure of Jordan
Canonical Form (JCF) in terms of the Segre characteristics is called a 
{\em bundle} that is proved to be a manifold \cite{Arn71,Gib76} 
through differential geometry.
Bundles can be established as manifolds using the geometric modeling approach
and Theorem~\ref{t:cam} but the complete proof is beyond the scope of this 
paper.
We illustrate the geometric modeling of a bundle using a specific 
JCF structure here.
Let 
\[
\Pi = \left.\left\{ X\,J_n(\la) \,X^{-1} \,\right|\, \la\,\in\,\CC,
X\,\in\,\CC^{n\times n} \,\,\mbox{is invertible}\right\}
\]
where $J_n(\la)$ ~denotes the $n\times n$ elementary Jordan block with the 
eigenvalue ~$\la$.
Namely $\Pi$ is the collection of all $n\times n$ matrices with a single 
eigenvalue in a single Jordan block.
The JCF problem with respect to this Jordan structure can be modeled as 
follows. 
At any $A\in\Pi$, pick a random vector $\bdc\in\CC^n$.
For almost all such $\bdc$, there is a unique invertible matrix 
$X\in\CC^{n\times n}$ whose columns are eigenvector and generalized 
eigenvectors such that $A\,X=X\,J_n(\la_*)$ along with the auxiliary equation
$\bdc^\h\,X=[1,0,\cdots,0]$.
\[ \mbox{Solve\,\, }\bdg(A,\la,Z) = (\bdo, O) \mbox{\,\, for \,\,}
(\la,Z)\,\in\,\CC\times\CC^{n\times n}
\]
with the holomorphic mapping from
\[ \mapform{\bdg}{\Omega\subset\CC^{n\times n}\times\CC\times\CC^{n\times n}}{
\CC^{1\times n}\times\CC^{n\times n}}{(G,\la,Z)}{
\big(\bdc^\h Z - [1,0,\cdots,0], G\,Z-Z\,J_n(\la)\big)}
\]
where $\Omega$ is a neighborhood of $(A,\la_*,X)$ in which all 
$(G,\la,Z)$ has an invertible $Z$ and nonzero dot-product between 
$\bdc$ and the lone eigenvector of $G$.
The fact that the subset $\Pi$ is a manifold can be established by verifying 
the four conditions in Theorem~\ref{t:cam} on $\bdg$ using common techniques 
in linear algebra, and
\[ \codim(\Pi) \,=\, \dm(\CC^{1\times n}\times\CC^{n\times n})
-\dm(\CC\times\CC^{n\times n}) \,=\, n-1.
\]

\vspace{-4mm}
\section{The least squares problem}

\vspace{-4mm}
As elaborated in \S\ref{s:mod}, an algebraic problems can be modeled as a 
zero-finding problem in the form of $\bdf(\bdu,\bdv)\,=\,\bdo$ for the 
variable $\bdv$ at a certain fixed parameter $\bdu$, and the equation is 
often overdetermined.
In practical computation, the parameter $\bdu$ is expected to be represented 
via empirical data $\tilde\bdu$ at which the exact solution $\bdv$ generally 
does not exist for the perturbed equation $\bdf(\tilde\bdu,\bdv)=\bdo$.
The resulting model becomes a least squares problem.

Let $\cV$ and $\cW$ be normed vector spaces isometrically isomorphic to 
$\CC^n$ and $\CC^m$ respectively with $m>n$. 
Let $\bdx\mapsto\bdf(\bdx)$ be a mapping from an open subset $\Omega$ of $\cV$ 
to $\cW$. 
Since $\bdf(\Omega)$ is of dimension at most $n$ in $\cW$ with $\dm(\cW)=m>n$, 
conventional solutions to the equation $\bdf(\bdx)=\bdb$ do not exist in 
general.
Instead, we seek a {\em least squares solution} $\bdx_*\in\Lambda$ of 
~$\bdf(\bdx)=\bdb$ such that 
\[ \big\|\bdf(\bdx_*)-\bdb\big\|^2 \,=\,
\min_{\bdx\in\Lambda}\big\|\bdf(\bdx)-\bdb\big\|^2
\]
where $\Lambda\subset\Omega$ is an open neighborhood of $\bdx_*$.
In other words, we seek $\bdx_*$ so that $\bdf(\bdx_*)$ is the projection of 
$\bdb$ to the surface $\bdf(\Omega)$, minimizing the distance from 
$\bdb$ to $\bdf(\Omega)$.
Further assume $\cV$ and $\cW$ are isometrically isomorphic to $\CC^n$ 
and $\CC^m$ respectively so that $\bdf_\bdx(\bdz)^\h$ and 
$\bdf_\bdx(\bdz)^\dagger$ are well defined.
Then a least squares solution is a {\em critical point} for the equation
$\bdf(\bdx)=\bdb$, namely (c.f. \cite{zeng-mr-05})
\begin{equation}\label{cp}   
\bdf_\bdx(\bdx_*)^\h\,\big(\bdf(\bdx_*)-\bdb\big) = \bdo.
\end{equation}
The {\em Gauss-Newton iteration}%
\footnote{A general purpose MATLAB module {\tt GaussNewton} is implemented 
in the package {\sc NAClab} \cite{naclab} with an intuitive interface
\cite{solve}.}
\begin{equation}\label{gnit}
	\bdx_{k+1} = \bdx_k - \bdf_\bdx(\bdx_k)^\dagger 
	\big(\bdf(\bdx_k)-\bdb\big) 
	\mbox{\,\,for\,\,} k = 0, 1, \ldots
\end{equation}
is effective in finding the least squares solution of $\bdf(\bdx)=\bdb$ 
and is locally convergent.
The following lemma provides detailed convergence conditions 
in Kantorovich style. 

\vspace{-4mm}
\begin{lemma}{\em \cite{ZengAIF}}  \label{gnlem}
Let $\cV$ and $\cW$ be finite-dimensional normed vector spaces 
isometrically isomorphic to $\CC^n$ and $\CC^m$ respectively.
Assume $\bdx\mapsto\bdf(\bdx)$ is a holomorphic mapping
from an open domain $\Omega\subset \cV$ to $\cW$ with a
critical point $\bdx_* \in \Omega$ of the system 
$\bdf(\bdx)=\bdb$ and $\bdf_\bdx(\bdx_*)$ is injective.
Then there is an open neighborhood $\Lambda\subset\Omega$ of $\bdx_*$ 
along with constants $\zeta, \,\gamma > 0$ such that 
\begin{equation} \label{Jc1}
\big\| \bdf_\bdx(\bdz)^\dagger \big\| \le \zeta 
\mbox{~and~}
\big\| \bdf(\bdz) - \bdf(\tilde\bdz) - 
\bdf_\bdx(\tilde\bdz)(\bdz-\tilde\bdz) \big\| 
\le \gamma\,\big\|\bdz-\tilde\bdz\big\|^2 
\end{equation}
for all $\bdz,\,\tilde\bdz\in \Lambda$.
Further assume $\|\bdf(\bdx_*)-\bdb\|$ is small so that
\begin{equation} \label{resloc}
\big\|\big( \bdf_\bdx(\bdz)^\dagger -\bdf_\bdx(\bdx_*)^\dagger \big) 
\big(\bdf(\bdx_*)-\bdb \big) \big\| \le
\sg\, \big\|\bdz-\bdx_*\big\|
\end{equation}
for a constant $\sg < 1$ at every $\bdz \in \Lambda$.
Then for all $\bdx_0 \in \Lambda$ satisfying 
\[
\big\|\bdx_0-\bdx_*\big\| < \mbox{$\frac{1-\sg}{\zeta \gamma}$}
\,\,\,\mbox{and}\,\,\, 
\{\bdx \in \cV \,|\, \|\bdx-\bdx_*\| < \|\bdx_0 - \bdx_*\|\}\subset
\Lambda, 
\]
the iteration {\em (\ref{gnit})} is well defined in $\Lambda$, converges to 
$\bdx_*$, and satisfies
\[ \big\| \bdx_{k+1} - \bdx_* \big\| \leq 
\big(\sg \, + \zeta \gamma \, \big\| \bdx_k - \bdx_* \big\|\big)\,
\big\| \bdx_k - \bdx_* \big\| 
\]
for $k=0,1,\ldots$ with $\sg + \zeta \gamma\, \big\|\bdx_0-\bdx_*\big\| < 1$.
\end{lemma}

\vspace{-6mm}
\section{Tubular Neighborhood Theorem}
\label{s:tnt}

\vspace{-4mm}
The very reason we need to establish manifolds in regularizing ill-posed 
algebraic problems lies in one of the fundamental theorems in differential 
geometry:
A smooth manifold is contained in an open {\em tubular neighborhood}\,
in which every point can be uniquely projected onto the manifold following 
a normal line and the projection mapping possesses certain desired properties.
The concept of tubular neighborhood is also regarded as ``one of the
most useful notions in the theory of differential manifolds'' \cite{Dieud}.
Standard versions of the tubular neighborhood theorem for {\em real} smooth 
manifolds can be found in textbooks of differential geometry (see e.g. 
\cite{BurGid}).
Those versions are in abstract forms for general purposes and
do not appear to be applicable to our geometric modles involving 
{\em complex} analytic manifolds. 
For the applications in regularization of ill-posed algebraic problems,
the projection to the manifold does not need to be holomorphic and it 
suffices to be Lipschitz continuous with the Lipschitz constant serving as a 
condition number measuring the sensitivity of the underlying problem.

\vspace{-4mm}
\begin{lemma}\label{l:lstnt}
Let $\cU$, $\cV$ and $\cW$ be normed vector spaces over $\CC$ that are 
isometrically isomorphic to $\CC^l$, $\CC^m$ and $\CC^n$ respectively 
with $m\le n\le l+m$.
Assume $\Pi$ is a complex analytic manifold in $\cU$ and, for every 
$\bdu_0\,\in\,\Pi$, there is a holomorphic mapping
$(\bdu,\bdv)\,\mapsto\,\bdf(\bdu,\bdv)$ from an open domain $\Omega\subset
\cU\times\cV$ to $\cW$ satisfying the conditions {\em (i)-(iv)} in
Theorem~\ref{t:cam}.
Then the following assertions hold:

{\em (i)} There are open neighborhoods $\Psi$ of $\bdu_0$ in $\cU$ and 
$\Phi$ of $\bdv_0$ in $\cV$ along with a mapping 
$\pi : \Psi\subset\cU \rightarrow \cV$ 
whose image $\tilde\bdv = \pi(\tilde\bdu) \in \Pi$ is the unique least 
squares solution to the equation $\bdf(\tilde\bdu,\bdv) = \bdo$ in $\Phi$ 
at every $\tilde\bdu \in \Psi$.
Furthermore, for every open neighborhood $\check\Phi\subset\Phi$ of 
$\bdv_0$ in $\cV$, there is an open neighborhood $\check\Psi\subset
\Psi$ of $\bdu_0$ such that $\pi(\check\Psi)\subset\check\Phi$.

{\em (ii)} The mapping $\pi$ is locally Lipschitz continuous in $\Psi$.

{\em (iii)} From every $\tilde\bdu\in\Psi$ serving as empirical 
data for $\bdu_0$, the least squares solution $\pi(\tilde\bdu)=
\tilde\bdv$ satsifies
\begin{equation}\label{Dlbdv}
  \|\tilde\bdv-\bdv_0\| \le \|\bdf_\bdv(\bdu_0,\bdv_0)^\dagger\|\,
\|\bdf_\bdu(\bdu_0,\bdv_0)\|\,\|\tilde\bdu-\bdu_0\| + 
o(\|\tilde\bdu-\bdu_0\|)
\end{equation}
\end{lemma}

\vspace{-4mm}
{\em Proof.}
Using the notations in the proof of Theorem~\ref{t:cam}, there exists a bounded 
open neighborhood $\Sigma$ of $\bdv_0$ in $\cV$ such that the subset
$\{\hat\bdu_0\}\times(\{\check\bdu_0\}\times\overline{\Sigma})\subset
\Lambda\times\Dl$. 
For any $r>0$ and the subset $\Phi_r:=\{\bdv\in\Sigma|\|\bdv-\bdv_0\|<r\}$,
we claim there is an $s>0$ such that, at every 
$\tilde\bdu\in \Psi_s :=\{\bdu\in\cU|\|\bdu-\bdu_0\|<r\}$, the minimum 
$\min_{\bdv\in\overline{\Phi_r}} \|\bdf(\tilde\bdu,\bdv)\|$ occurs at a 
certain $\tilde\bdv\in\Phi_r$ that is a least squares solution of 
$\bdf(\tilde\bdu,\bdv)=\bdo$.
Assume otherwise. 
Then there is a sequence $\{\bdu_j\}_{j=1}^\infty$ converging to $\bdu_0$ 
such that $\min_{\bdv\in\overline{\Phi_r}} \|
\bdf(\bdu_j,\bdv)\|=\|\bdf(\bdu_j,\bdv_j)\|$ at $\bdv_j\in
\overline{\Phi_r}\setminus\Phi_r$ for every $j=1,2,\ldots$.
Since $\overline{\Phi_r}\setminus\Phi_r$ is compact, we can assume
$\bdv_j$ converges to a certain $\check\bdv$.
Thus 
\[ \|\bdf(\bdu_0,\check\bdv)\| = 
\lim_{j\rightarrow\infty}\|\bdf(\bdu_j,\bdv_j)\| \le 
\lim_{j\rightarrow\infty}\|\bdf(\bdu_j,\bdv_0)\| = 0,
\]
implying $\check\bdv\,=\,\bdv_0$ that contradicts to 
$\check\bdv\,\in\,\overline{\Phi_r}\setminus\Phi_r$.

We can assume $r_1 > 0$ is sufficiently small so that, 
for every $\bdv_1, \bdv_2 \in \Phi_{r_1}$ and $\bdu \in \Psi_{s_1}$, 
there exist constants $\zeta, \gamma > 0$ such that
\begin{align*}
\|\bdf(\bdu,\bdv_2)-\bdf(\bdu,\bdv_1)-\bdf_\bdv(\bdu,\bdv_1)\,(\bdv_2-\bdv_1)\|
&< \gamma\,\|\bdv_2-\bdv_1\|^2 \\
\big\|\big(\bdf_\bdv(\bdu,\bdv_2)^\dagger-\bdf_\bdv(\bdu,\bdv_1)^\dagger\big)\,
\bdf(\bdu,\bdv_1)\big\| &<\mbox{$\frac{1}{2}$}\,\|\bdv_2-\bdv_1\| \\
\|\bdf_\bdv(\bdu,\bdv_1)^\dagger\| < \zeta, \,\,\,\,\,
\|\bdv_2-\bdv_1\| &< \mbox{$\frac{1}{2\,\zeta\,\gamma}$}.
\end{align*}
Let $r_2\,=\,\frac{1}{3}\,r_1$, $\Phi\,=\,\Phi_{r_2}$ and 
$\Psi\,=\,\Psi_{s_1}\cap\Psi_{s_2}$.
For every $\hat\bdu\,\in\,\Psi$, the minimum 
$\min_{\bdv\in\overline{\Phi}}\|\bdf(\hat\bdu,\bdv)\|$ is attainable at a 
certain $\hat\bdv\,\in\,\Phi$ and, for any initial iterate $\bdv_1\,\in\,\Phi$,
we have $\|\bdv_1-\hat\bdv\|\,<\,\frac{1}{2\,\zeta\,\gamma}\,=\,
(1-\frac{1}{2})\frac{1}{\zeta\,\gamma}$ and the set
$\Omega\,=\,\{\bdv\in\cV | \|\bdv-\hat\bdv\|\,<\,\|\bdv_1-\hat\bdv\|\}$ 
is a subset of $\Phi_{r_1}$ since, for every $\bdv\,\in\,\Omega$, we have
\begin{align*}
\|\bdv-\bdv_0\| &\le \|\bdv-\hat\bdv\|+\|\hat\bdv-\bdv_0\| 
< \|\bdv-\hat\bdv\|+r_2 < \|\bdv_1-\hat\bdv\|+r_2 \\
&\le \|\bdv_1-\bdv_0\|+\|\bdv_0-\hat\bdv\|+r_2 
< r_2+r_2+r_2 = r_1 
\end{align*}
By Lemma~\ref{gnlem}, for every initial iterate $\bdv_1\,\in\,\Phi$, the
Gauss-Newton iteration on the equation $\bdf(\hat\bdu,\bdv)\,=\,\bdo$ 
converges to $\hat\bdv$.
This local minimum is unique in $\Phi$ because, assuming there is another
minimum point $\check\bdv\,\in\,\Phi$ of $\|\bdf(\hat\bdu,\bdv)\|$,
the Gauss-Newton iteration converges to $\hat\bdv$ from the initial iterate
$\check\bdv$. 
On the other hand, the Gauss-Newton iteration from the local minimum point
$\check\bdv$ stays at $\check\bdv$, implying $\check\bdv\,=\,\hat\bdv$ 
and thus the existence of the mapping $\pi$.
Given any open subset $\check\Phi$ of $\Phi$, there is an open subset
$\check\Psi$ of $\Psi$ for the same reason that $\Psi_s$ exists such that the 
minimum $\min_{\bdv\in\overline{\check\Phi}}\|\bdf(\tilde\bdu, \bdv)\|$ is 
attainable at a certain $\tilde\bdv\,\in\,\check\Phi$ for every fixed 
$\tilde\bdu\,\in\,\check\Psi$.
This $\tilde\bdv$ is unique in $\Phi$ since $\tilde\bdu\,\in\,\Psi$, 
and thus $\tilde\bdv$ is unique in $\check\Phi$, implying $\tilde\bdv
\,=\,\pi(\tilde\bdu)$ so that $\pi(\check\Psi)\,\subset\,\check\Phi$.

On the Lipschitz continuity the mapping $\pi$, let $\tilde\bdu,\,\hat\bdu
\,\in\,\Psi$ with $\pi(\tilde\bdu)\,=\,\tilde\bdv$ and $\pi(\hat\bdu)\,=\,
\hat\bdv$.
The one-step Gauss-Newton iteration 
$\bdv_1\,=\,\tilde\bdv-\bdf_\bdv(\hat\bdu,\tilde\bdv)^\dagger\,
\bdf(\hat\bdu,\tilde\bdv)$ from $\tilde\bdv$ on the equation 
$\bdf(\hat\bdu,\bdv)\,=\,\bdo$ toward $\hat\bdv$ yields the iniquality
$\|\bdv_1-\hat\bdv\|\,\le\,\mu\,\|\tilde\bdv-\hat\bdv\|$ with 
$\mu\,<\,1$ by Lemma~\ref{gnlem}.
Thus
\[ \|\tilde\bdv-\hat\bdv\| \le \|\hat\bdv-\bdv_1\|+\|\bdv_1-\tilde\bdv\|
\le \mu\,\|\tilde\bdv-\hat\bdv\|+\|\bdv_1-\tilde\bdv\|.
\]
Using the identity $\bdf_\bdv(\tilde\bdu,\tilde\bdv)^\dagger\,
\bdf(\tilde\bdu,\tilde\bdv)\,=\,\bdo$, the Lipschitz continuity of 
$\bdf$ and $\bdf_\bdv$ along with
\begin{align*}
&\|\bdf_\bdv(\hat\bdu,\tilde\bdv)^\dagger- 
\bdf_\bdv(\tilde\bdu,\tilde\bdv)^\dagger\| \\
&\le 3\,\|\bdf_\bdv(\tilde\bdu,\tilde\bdv)^\dagger\|^2\,
\|\bdf_\bdv(\hat\bdu,\tilde\bdv)- \bdf_\bdv(\tilde\bdu,\tilde\bdv)\|
\tag{c.f. \cite[Theorem 3.4]{stew77}}  \\
&\le 3\,\|\bdf_\bdv(\tilde\bdu,\tilde\bdv)^\dagger\|^2\,
\|(\bdf_\bdv)_\bdu(\tilde\bdu,\tilde\bdv)\|\,\|\hat\bdu-\tilde\bdu\|
+ O(\|\hat\bdu-\tilde\bdu\|^2)
\end{align*}
for sufficiently small $\|\hat\bdu-\tilde\bdu\|$ where
$(\bdf_\bdv)_\bdu(\tilde\bdu,\tilde\bdv)$ is the Jacobian of the 
holomorphic mapping $\bdu\,\mapsto\,\bdf_\bdv(\bdu,\tilde\bdv)$ at 
$\tilde\bdu$, we have
\begin{eqnarray}
\lefteqn{\|\tilde\bdv-\hat\bdv\|
\le \mbox{$\frac{1}{1-\mu}$}\,\|\bdv_1-\tilde\bdv\|} \nonumber
\\
&=&\mbox{$\frac{1}{1-\mu}$}\,\|\bdf_\bdv(\hat\bdu,\tilde\bdv)^\dagger\,
\bdf(\hat\bdu,\tilde\bdv)-\bdf_\bdv(\tilde\bdu,\tilde\bdv)^\dagger\,
\bdf(\tilde\bdu,\tilde\bdv)\|  \nonumber \\
&\le&\mbox{$\frac{1}{1-\mu}$}\,
\big(\|\bdf_\bdv(\tilde\bdu,\tilde\bdv)^\dagger\|\,
\|\bdf(\hat\bdu,\tilde\bdv)-\bdf(\tilde\bdu,\tilde\bdv)\| 
+ \|\bdf_\bdv(\hat\bdu,\tilde\bdv)^\dagger-
\bdf_\bdv(\tilde\bdu,\tilde\bdv)^\dagger\|\,
\|\bdf(\hat\bdu,\tilde\bdv)\|\big) \nonumber \\
&\le& \mbox{$\frac{\|\bdf_\bdv(\tilde\bdu,\tilde\bdv)^\dagger\|
}{1-\mu}$}\, 
\big(
\|\bdf_\bdu(\tilde\bdu,\tilde\bdv)\|
+3\,\|\bdf_\bdv(\tilde\bdu,\tilde\bdv)^\dagger\|\,
\|(\bdf_\bdv)_\bdu(\tilde\bdu,\tilde\bdv)\|\,
\|\bdf(\tilde\bdu,\tilde\bdv\|\big)\,
\|\hat\bdu-\tilde\bdu\|\, \nonumber \\ & & 
 +O(\|\hat\bdu-\tilde\bdu\|^2).
\label{hutu}
\end{eqnarray}
As a result, there is a constant $\theta>0$ such that
$\|\tilde\bdv-\hat\bdv\|\,\le\,\theta\,\|\tilde\bdu-\hat\bdu\|$
when $\|\tilde\bdu-\hat\bdu\|$ is sufficiently small, leading to the 
assertion (ii). 
Set
\[ (\tilde\bdu,\tilde\bdv)\,=\,(\bdu_0,\bdv_0) ~~~~\mbox{and}~~~
(\hat\bdu,\hat\bdv)\,=\, (\bdu_0+\Dl\bdu,\bdv_0+\Dl\bdv)
\]
in (\ref{hutu}) and apply
$\bdf(\bdu_0,\bdv_0)\,=\,\bdo$ and $\mu\,=\,O(\|\hat\bdu-\tilde\bdu\|)$.
The inequality (\ref{Dlbdv}) holds.
\qed

Based on Lemma~\ref{l:lstnt}, the following Theorem~\ref{t:tnt} is a version 
of the Tubular Neighborhood Theorem for manifolds in normed 
vector spaces isometrically isomorphic to $\CC^n$'s.
It is specifically tailored for the application of solving ill-posed 
algebraic problems from empirical data.
There appears to be no such a version in the literature of differential
geometry since some analytic structures can not be preserved in
the tubular neighborhood and not needed in our application.
We provide a proof based on the Gauss-Newton iteration and Lemma~\ref{l:lstnt}.

\vspace{-4mm}
\begin{theorem}[Tubular Neighborhood Theorem] \label{t:tnt}
Let $\Pi$ be a complex analytic \linebreak
manifold in a vector space $\cU$ that is 
isometrically isomorphic to $\CC^n$.
There is a {\em tubular neighborhood}, namely an open subset 
$\Omega\,\supset\,\Pi$ of $\cU$ such that every $\bdb \in \Omega$ has a 
unique projection $\bdx_\bdb \in \Pi$ of minimum distance to $\bdb$, that is
\begin{equation}
\big\|\bdx_\bdb-\bdb\big\| = 
\inf_{\bdx \in \Pi} \big\|\bdx-\bdb\big\| =: \dist{\bdb,\,\Pi}.
\end{equation}
Furthermore, the projection $\bdb\,\mapsto\, \bdx_\bdb$ from $\Omega$ to $\Pi$ 
is locally Lipschitz continuous.
\end{theorem}

\vspace{-4mm}
{\bf Proof}.
Let $\bdu_0$ be any particular point in $\Pi$.
Since $\Pi$ is a complex analytic manifold in $\cU$, there is an open 
neighborhood $\cM$ of $\bdu_0$ in $\cU$, an open subset $\cN$ of $\CC^m$ and 
a holomorphic mapping
$\bdv\,\mapsto\,\phi(\bdv)$ from $\cN\subset\CC^m$ onto $\cM\cap\Pi$ 
with a holomorphic inverse $\phi^{-1}$ from $\cM\cap\Pi$ onto $\cN$.
Let the holomorphic mapping $\bdf\,:\,(\bdu,\bdv)\,\mapsto\,
\phi(\bdv)-\bdu$ from $\cM\times\cN\,\subset\,\cU\times\CC^m$ to $\cU$.
Then $\bdf$ satisfies all the conditions of Lemma~\ref{l:lstnt}.
As a result, there is an open neighborhood $\Psi\subset\cM$ of $\bdu_0$ 
in $\cU$ such that, for every $\hat\bdu\,\in\,\Psi$, there exists a unique
least squares solution $\bdv\,=\,\hat\bdv$ for the equation 
$\bdf(\hat\bdu,\bdv)\,=\,\bdo$ so that 
\[ \|\bdf(\hat\bdu,\hat\bdv)\| \,=\min_{\bdv\in\Phi}\|\bdf(\hat\bdu,\bdv)\|
\,=\, \min_{\bdv\in\Phi}\|\phi(\bdv)-\hat\bdu\|.
\]
We can assume the neighborhood $\Psi$ is sufficiently small so that any
$\hat\bdu\,\in\,\Psi$ satisfies the inequality $\|\hat\bdu-\bdu\|\,>\,
\|\phi(\hat\bdv)-\hat\bdu\|$ for all $\bdu\,\in\,\Pi\setminus\phi(\Phi)$,
implying the local minimum $\|\phi(\hat\bdv)-\hat\bdu\|\,=\,
\min_{\bdu\in\Pi}\|\bdu-\hat\bdu\|$ is the global minimum.
\qed

From computational point of view, the desired solution $\hat\bdv\,\in\,\cV$ 
at a data point $\hat\bdu\,\in\,\cU$ may be modeled as the zero of a 
holomorphic mapping $\bdv\,\mapsto\,\bdf(\hat\bdu,\bdv)$ with $\hat\bdu$ in a 
structure-preserving manifold $\Pi$ in $\cU$. 
When $\hat\bdu$ is not known exactly but represented by its empirical data in
$\tilde\bdu\,\approx\,\hat\bdu$, the Tubular Neighborhood Theorem ensures 
the projection $\check\bdu$ of $\tilde\bdu$ to $\Pi$ uniquely exists, 
enjoys Lipschitz continuity and is independent of choices of the mapping
$\bdf$ in the model.
As a result, the solution $\check\bdv$ at the parameter value $\check\bdu$ 
can be defined as the {\em regularized solution} at $\tilde\bdu$. 
From Lemma~\ref{l:lstnt}, the regularized solution $\check\bdv$ can be 
accurately approximated by the least squares solution $\tilde\bdv$ of the 
equation $\bdf(\tilde\bdu,\bdv)\,=\,\bdo$ as long as $\bdf$ is properly 
constructed following the Geometric Modeling Theorem.

\vspace{-4mm}
\section{The geometric regularization: Concluding remarks and examples}

\vspace{-4mm}
As a notion attributed to Hadamard, a mathematical model is a 
{\em well-posed problem} if its solution satisfies existence, uniqueness and 
Lipschitz continuity with respect to data perturbations. 
Those problems may also be loosely referred to as being {\em regular}.
Otherwise, the problem is ill-posed or often called {\em singular}.
In general, singular problems are difficult to solve accurately from 
empirical data and require some form of regularization.

Algebraic problems such as polynomial GCD/factorizations,
matrix rank/kernels and matrix Jordan Canonical Forms (JCF) are not all 
singular.
For each problem, the data space is partitioned by manifolds and, on every 
manifold, the solutions maintains a specific algebraic structure.
Data associated with regular problems are open dense in the data
space, forming a manifold of codimension zero.
A problems is singular when the data point lies on a manifold 
of a positive codimension.
Due to the dimension deficit, a perturbation generically pushes the data away 
from the native manifold and the alters the structure of the solution, 
implying the solution is highly sensitive to {\em arbitrary} data 
perturbations.

However, the solutions of those singular problems are locally Lipschitz 
continuous if the data are constrained on a structure-preserving manifold.
By the Geometric Modeling Theorem, algebraic problem on any such manifold
$\Pi$ can be modeled as a zero finding problem ~$\bdf(\bdu,\bdv)=\bdo$
for the variable $\bdv$ at $\bdu\in\Pi$.
If an underlying data point ~$\hat\bdu$ is known with limited accuracy
through empirical data ~$\tilde\bdu$, Lemma~\ref{l:lstnt} and the 
Tubular Neighborhood Theorem (Theorem~\ref{t:tnt}) ensure that
solving for the least squares solution ~$\tilde\bdv$ of the equation
$\bdf(\tilde\bdu,\bdv)=\bdo$ is a well-posed problem. 
In other words, a singular algebraic problem can be regularized if it can be 
properly modeled following the Geometric Modeling Theorem assuming the 
structure of the solution is known.

Detailed elaboration on the identification of the solution structure is 
beyond the scope of this paper. 
In a nutshell, we can quantify the {\em singularity} of each data point as 
the codimension of the manifold on which the data point resides. 
The structure-preserving manifolds are entangled to form a strata in which
every manifold is embedded in the closures of some manifolds of lower 
codimensions.
In other words, such an algebraic problem is highly sensitive but the 
sensitivity is {\em directional} such that sufficiently small perturbations 
can only reduce the singularity and never increase it.

At an underlying data point ~$\hat\bdu$ on a structure-preserving manifold
$\Pi$, the given empirical data point ~$\tilde\bdu$ is a small perturbation 
from ~$\hat\bdu$.
Assuming the perturbation is sufficiently small so that $\tilde\bdu$ stays
in the tubular neighborhood, the underlying manifold $\Pi$ is 
of the highest singularity (codimension) among all the manifolds intersect
a small neighborhood of ~$\tilde\bdu$. 
Identification of the solution structure becomes a discrete optimization 
problem in maximizing the codimension (singularity) of the manifolds within
an error tolerance of the empirical data point $\tilde\bdu$.
Consequently, a natural strategy for computing the regularized solution
at an empirical data parameter $\tilde\bdu$ is a two-staged process in 
either symbolic, numerical or hybrid computation:

\vspace{-4mm}
\begin{itemize}\parskip0mm
\item[]\hspace{-5mm}{\bf Stage I}. Within an error tolerance of the data 
$\tilde\bdu$, 
find the nearby structure-preserving manifold of the highest singularity.
\item[] \hspace{-5mm}{\bf Stage II}.
Solve the equation $\bdf(\tilde\bdu,\bdv)=\bdo$ that is properly formulated
based on the Geometric Modeling Theorem for its least squares solution 
$\bdv=\tilde\bdv$.
\end{itemize}

\vspace{-4mm}
The least squares solution ~$\tilde\bdv$ is a {\em regularized solution} 
at the empirical data $\tilde\bdu$ {\em within the error tolerance}.
It is {\em not} a solution at ~$\tilde\bdu$ in conventioanl sense but 
accurately approximates the exact solution at the underlying data $\hat\bdu$
with an error bound proportional to the data error ~$\|\tilde\bdu-\hat\bdu\|$.

This regularization strategy has been applied to many singular algebraic 
problems such as computing multiple roots of univariate polynomials
\cite{ZengAIF}, approximate polynomial GCD \cite{ZengGCD,zeng-dayton}, 
factorization of multivariate polynomials \cite{WuZeng} from empirical data. 
The resulting algorithms are implemented in the MATLAB package {\sc NAClab}
\cite{naclab} including a preliminary module for computing the Jordan Canonical 
Form from possibly perturbed matrices.
We illustrate the strategy the following examples.

\vspace{-4mm}
\begin{example}\label{e}\em
Let the polynomial pair $(p,q)\in\mP_{13}\times\mP_{11}$ where
\begin{align*}
\tilde{p} &\,=\, 
1 - 0.333 x + 0.667 x^3 + x^{10} - 0.333 x^{11} + 0.666 x^{13} \\
\tilde{q} &\,=\, 1.429 + 3.571 x + 1.429 x^{10} + 3.571 x^{11}
\end{align*}
serving as empirical data of the pair $(p,q)$ that equals
\[
\left( 1 - \mbox{$\frac{1}{3}$} x + \mbox{$\frac{2}{3}$} x^3 + x^{10} 
- \mbox{$\frac{1}{3}$} x^{11} + \mbox{$\frac{2}{3}$} x^{13}, \,
\mbox{$\frac{10}{7}$} + \mbox{$\frac{25}{7}$} x + 
\mbox{$\frac{10}{7}$} x^{10} - \mbox{$\frac{25}{7}$} x^{11}\right).
\]
In exact sense, we have $\GCD(p,q)=1+x^{10}$ but $\GCD(\tilde{p},\tilde{q})=1$
that are far apart due to the singularity of the GCD even though the data error is about $10^{-3}$.
The computing objective is to find an approximate GCD $\approx 1+x^{10}$ from
the empirical data $(\tilde{p},\tilde{q})$ by calculating a regularized GCD
within the data error bound $10^{-3}$.

At Stage I, the GCD degree (i.e. structure) is identified by computing 
the numerical nullity of the Sylvester matrix $S(\tilde{p},\tilde{q})$ within
the error tolerance $10^{-3}$. 
This numerical nullity is identical to the degree 10 of $\GCD(p,q)$.
Therefore, the native GCD manifold is $\cP_{13,11}^{10}$.
Further more, initial approximation $(u_0,v_0,w_0)$ of the GCD and cofactors
can be obtained by solving two linear systems (c.f. \cite{ZengGCD}).

At Stage II, we formulate the geometric model by constructing the mapping
$\psi$ as in (\ref{gcdmap}) for $k=10$, $m=13$, $n=11$ and solve the equation
$\psi(u,v,w,\tilde{p},\tilde{q}) \, =\, (0,0,0)$ for the least squares 
solution $(u,v,w)\in\mP_2\times\mP_{11}\times\mP_9$ with $u\sim 
1+0.9998 x^{10}$ with an accuracy in the order of the data error bound.
The regularized GCD computation is implemented in {\sc NAClab} so that
the above computation can be carried out in simple MATLAB sequence:

\scriptsize
\begin{verbatim}
>> p = '1-.333*x+0.667*x^3+x^10-0.333*x^11 + 0.666*x^13'; % enter polynomial p
>> q = '-1.429 - 3.571*x - 1.429*x^10 - 3.571*x^11';      % enter polynomial q
>> pgcd(p, q, 0.001) % regularized GCD of p and q within error tolerance 0.001
ans =
    '-1.24459473398662 - 1.24432753501985*x^10'
\end{verbatim}
\normalsize

\noindent The result is a multiple of $1+0.9998 x^{10}$.
\end{example}

\vspace{-4mm}
Computing Jordan Canonical Forms of matrices from empirical data is known 
to be a tremendous challenge.  
We conclude this paper with two examples of the module {\tt RegularizedJCF}
in {\sc NAClab} based on the geometric modeling elaborated in this paper. 

\vspace{-4mm}
\begin{example} \em
~There are applications where empirical data may even be
preferred over exact ones.
Consider the matrix
\[ A(r,s,t) =
 \mbox{\scriptsize $\left[\begin{array}{rrrrrr} 
\mbox{\em 2r-2s+t} & \mbox{\em 1-s+t} & \mbox{\em r-3-3s+2t} & 
\mbox{\em r-1-2s+t} & \mbox{\em -1} & \mbox{\em -1-s+t} \\
\mbox{\em r+4+s-2t} & \mbox{\em 3r+2-2t}& \mbox{\em 2r+10+2s-4t}& 
\mbox{\em r+5+s-2t}& \mbox{\em -r+s}& \mbox{\em r+3+s-2t} \\
\mbox{\em 1+4r-3s-t}& \mbox{\em 1+3r-2s-t}& \mbox{\em 1+7r-4s-2t}& 
\mbox{\em 1+4r-3s-t}& \mbox{\em -r-1+s}& \mbox{\em 2r-s-t} \\
\mbox{\em 7s-t-6r-2} &\! \mbox{\em 4s-t-3r-3}& \mbox{\em 10s-2t-8r+1}& 
\mbox{\em 7s-t-5r-1}& \mbox{\em 3+r-s}& \mbox{\em 3s-t-2r+1} \\
\mbox{\em r+3+2s-3t}\!&\! \mbox{\em 3r+1-3t}& \mbox{\em 2r+9+4s-6t}& 
\mbox{\em r+4+2s-3 t}& \mbox{\em 1-r+2s}& \mbox{\em r+3+2s-3t} \\
\mbox{\em -5r-4+5t}& \mbox{\em s+5t-6r-2}& \mbox{\em 10t-9r-10-s}& 
\mbox{\em -5r-5+5t}& \mbox{\em 2r-2s}& \mbox{\em 5t -3r-3-s}
\end{array}\right]$} 
\]
whose JCF is known to be ~$J_3(r)\oplus J_2(s) \oplus J_1(t)$.
When the parameter values $r$, $s$ and $t$ are exact, say 
{\tiny $\sqrt{k+\sqrt{k+\sqrt{k}}}$} for $k=2,3,5$, test on 
Maple 17 could not finish after hours of computation.
We can instead use approximation by rounding to, say 5 digits after decimal 
and find the regularized JCF within the error tolerance ~$10^{-4}$. 
The following is a demo of using {\tt RegularizedJCF} in {\sc NAClab}
that takes negligible elapsed time 0.03 second.

\vspace{-3mm}
{\scriptsize
\begin{verbatim}
>> A = [ 214636   149815  -231707   -81521  -100000   -50185 % enter matrix data
         269034   233854   738068   369034    31336   169034
         -75161   -43824     8509   -75161   -68664  -112488
        -061796  -255806   251061   234361   268664   112858
         119219  -143454   538438   219219   358830   119219
           5757   237093  -219823   -94243   -62673   270577]/100000;
>> [J,X] = RegularizedJCF(A,1e-4);   % call the software module
>> single(J)                         % display JCF in single precision
ans = 
   1.9615549   0.4031104           0           0           0           0
           0   1.9615549   3.7739313           0           0           0
           0           0   1.9615549           0           0           0
           0           0           0   2.2749500  -1.2751906           0
           0           0           0           0   2.2749500           0
           0           0           0           0           0   2.7730999
\end{verbatim}
}

\vspace{-4mm}
\noindent obtaining the exact JCF structure and eigenvalues of an accuracy
that is moderately proportional to that of the data.
\end{example}

\vspace{-4mm}
\begin{example}\em
~Godunov \cite[page 10]{godunov} uses the matrix
\[  
G ~~=~~ \left[ \mbox{\scriptsize $\begin{array}{rrrrrrr}
289 &  2064  &  336 &   128 &    80 &     32 &    16  \\
1152  &   30 &  1312 &   512 &   288 &    128 &    32 \\
-29 & -2000 &   756 &   384 &  1008 &    224 &    48 \\
512 &   128  &  640 &     0 &   640 &    512 &   128 \\
1053 &  2256 &  -504 &  -384 &  -756 &    800 &   208 \\
-287  &  -16 &  1712 &  -128 &  1968  &   -30 &  2032 \\
-2176 &  -287 & -1565 &  -512 &  -541 &  -1152 &  -289
\end{array}$} \right]  \]
to illustrate the difficulties in computing eigenvalues.
The eigenvalues $0,\pm 1, \pm 2, \pm 4$ of $G$ are simple but extremely 
sensitive with condition numbers arround $4\times 10^{12}$, implying $G$ is 
a small perturbation from a matrix on a singular bundle.
Trying an error tolerance, say $10^{-9}$, the module {\tt RegularizedJCF} in 
{\sc NAClab} finds the regularized JCF of $G$ within $10^{-9}$ as
a direct sum of of a $4\times 4$ and a $3\times 3$ elementary Jordan blocks
\[ \tilde{J}~~=~~ J_4(-2.121366210414752)\oplus J_3(2.828488280553040).
\]
This is the JCF of a nearby matrix $\hat{G}$ of singularity 5 with a 
relative distance 
\[ \mbox{$\frac{\|G-\hat{G}\|_{_F}}{\|G\|_{_F}}$} ~\approx~ 3.13\times 10^{-10}.
\]
The condition number of the JCF of $\hat{G}$ is much smaller at
$3.6\times 10^6$. 
There is another nearby bundle of even higher singularity.
Seting an error tolerance, say $0.005$, the regularized JCF of $G$ becomes a 
single $7\times 7$ elementary Jordan block $J_7(0.000000000001459)$.
with a moderate condition number 4268.5.
In 9-digit integer representation, there is a matrix $\tilde{G}=X\,J\,X^{-1}$ 
with an exact eigenvalue zero in a $7\times 7$ elementary Jordan block and a 
relative distance
\[ \mbox{$\frac{\|G-\tilde{G}\|_{_F}}{\|G\|_{_F}}$} ~\approx~ 5.3\times 10^{-7} 
\]
where
\begin{eqnarray*} 
X &~=~&
\left[\mbox{\tiny $\begin{array}{rrrrrrr}
\mbox{-500000494}\!\! &  \mbox{499231619}\!\! &  \mbox{475440501} &  \mbox{550430216} & \mbox{249476819} &  \mbox{344543097}  & \mbox{244097} \\
\mbox{244296}\!\!  & \mbox{39690036}\!\! &  \mbox{305346098} &  \mbox{418811015} &  \mbox{245808894} &  \mbox{229491349} & \mbox{499999756} \\
\mbox{499998993}\!\! & \mbox{-499191461}\!\! & \mbox{-425110651} &  \mbox{-12743120} &  \mbox{-32442421} &  \mbox{-25936159} &  \mbox{500000266} \\
\mbox{406}\!\!  &  \mbox{240}\!\! &  \mbox{-293020} &  \mbox{209706406} &  \mbox{260079479} &  \mbox{212082338}  &  \mbox{-33} \\
\mbox{-499998969}\!\! &  \mbox{499191323}\!\! &  \mbox{425111477} & \mbox{11515088} &  \mbox{494926230} &  \mbox{166312083} &  \mbox{500000239} \\
\mbox{500001425}\!\! & \mbox{-499232519}\!\! & \mbox{-475441212} & \mbox{-550432966} & \mbox{-249025474} &  \mbox{786194594}  & \mbox{244019} \\
\mbox{-244271}\!\! &  \mbox{-39689958}\!\! & \mbox{-305347057} & \mbox{-417583165} & \mbox{-708293582} & \mbox{-370419391} & \mbox{499999620}
\end{array}$}\right] \\
J & ~=~ &  \mbox{\tiny \,$\frac{1}{10^6}$}
\left[\mbox{\tiny $\begin{array}{rrrrrrr}
 0\!\! &  \mbox{-163589092} & 0 & 0 & 0 & 0 & 0 \\
 0\!\! & 0\!\! & \mbox{1279307109} &  0 & 0 & 0 & 0 \\
 0\!\! & 0\!\! & 0 &  \mbox{2151028721} & 0 & 0 & 0 \\
 0\!\! & 0\!\! & 0 & 0 &   \mbox{113025963} & 0 & 0 \\
 0\!\! & 0\!\! & 0 & 0 & 0 &  \mbox{2502078868} & 0 \\
 0\!\! & 0\!\! & 0 & 0 & 0 & 0 &\!\!  \mbox{3622414612} \\
 0\!\! & 0\!\! & 0 & 0 & 0 & 0 & 0
\end{array}$}\right]
\end{eqnarray*}
We now have a geometric interpretation on the sensitivity of $G$ in
an eigenproblem.
Let $\Pi_7$, $\Pi_{4\,3}$ and $\Pi_{1111111}$ ~be bundles corresponding
to Jordan structures $J_7(\la)$, $J_4(\la_1)\oplus J_3(\la_2)$ ~and
$J_1(\la_1)\oplus\cdots\oplus J_1(\la_7)$ respectively.
The bundle ~$\Pi_7$ is of the highest singularity (i.e. codimension) 6
among all bundels in $\CC^{7\times 7}$ and is embedded in the closure
of $\Pi_{4\,3}$ with a lower singularity 5 while both bundles are in the
closure of the open dense bundle $\Pi_{1111111}$ of singularity zero.
Although $G\in\Pi_{1111111}$ is regular, its eigenproblem is highly 
ill-conditioned  because $G$ is a tiny distance $10^{-10}$ from the bundle 
$\Pi_{4\,3}$ of singularity 5 and $10^{-7}$ from the most singular bundle
$\Pi_7$. 
With proper geometric modeling, the regularized JCF problem of $G$ is not
as ill-conditioned as the straightforward eigenproblem.
\end{example}

{\bf Acknowledgments.}
The author is indebted to his former colleague Marian Gidea for introducing 
the Tubular Neighborhood Theorem in a conversation leading to this work.

\end{document}